\documentclass[12pt]{article}

\usepackage[T1,T2A]{fontenc}
\usepackage[utf8]{inputenc}

\usepackage{amsfonts}
\usepackage{amssymb}
\usepackage{amsmath}
\usepackage{amsthm}

\def \le {\leqslant}
\def \ge {\geqslant}

\usepackage{tikz} 
\usetikzlibrary{calc,arrows,decorations.pathreplacing,fadings,3d,positioning}

\begin{document}

\begin{Huge}
\centerline{Über Approximationen  $n$ reeller Zahlen}
\end{Huge}
\vskip+0.5cm

\begin{Large}
\centerline{von Vassily Manturov und Nikolay Moshchevitin\footnote{
Dieser Autor wurde  durch RNF unterstützt   (Grant No.  19-11-00001).
}
}
\end{Large}
\vskip+2cm

Sei $\alpha$ eine reelle Zahl und
$$
\psi_\alpha (t) =\min_{q\in \mathbb{Z}_+, q\le t}  || q\alpha ||,\,\,\,\,\, \text{wohin}\,\,\,\,\,
||x || = \min_{a\in \mathbb{Z}}|x-a|
$$
die zugehörige Funktion des  Irrationalitätsmaßes. Im Jahr 2010 haben
 Kan und Moshchevitin \cite{KM} bewiesen, dass
   die Differenz  
 $\psi_\alpha (t) - \psi_\beta(t)$
 für zwei irrationalen  Zahlen $\alpha$ und $\beta$     
 unendlich oft   ihr Vorzeichen wechselt  wenn $ t\to \infty$.
 In letzter Zeit sind mehrere Artikel erschienen, die dieses Ergebnis auf unterschiedliche Weise verallgemeinern
 (siehe \cite{CKM, D, M_pal, M_arch, Shu, Sha}).
 Das Ziel dieser Arbeit ist eine neue
 Verallgemeinerung des Satzes von Kan und Moshchevitin  für $n\ge 2 $ reellen  Zahlen zu beweisen.
 \vskip+0.5cm
 
 {\bf 1. Über $\frak{k}$-Index.}
  \vskip+0.3cm
 
  Wir nennen   die Zahlen $\alpha, \beta\in \mathbb{R}$  {\it inkommensurabel},
  wenn die Ungleichung
  $$
  \psi_\alpha (t) \neq \psi_\beta(t)
  $$
  für alle hinreichend großen $t$ gilt. Es ist klar, dass wenn  $1, \alpha,\beta$ 
   linear unabhängig über $\mathbb{Q}$ sind, 
    die Zahlen $\alpha$ und $\beta$  inkommensurabel sind.
     Das Gegenteil ist   aber nicht der Fall.
   
   Sei  $\pmb{\alpha} = (\alpha_1,...,\alpha_n)$  ein $n$-Tupel  paarweise inkommensurabler Zahlen.
   Dann ist
   $
  \psi_{\alpha_i} (t) \neq \psi_{\alpha_j}(t)
  $
  für alle hinreichend großen $t$  und für alle $ i\neq j$. So gibt es für jede hinreichend große $t$ eine Permutation
  $$
 \pmb{\sigma} (t)  : \{1,2,3,...,n\} \mapsto  \{\sigma_1,\sigma_2,\sigma_3,...,\sigma_n\}
 $$
 mit
 $$
 \psi_{\alpha_{\sigma_1}} (t) >
 \psi_{\alpha_{\sigma_2}} (t) >
  \psi_{\alpha_{\sigma_3}} (t) >
  ... >
   \psi_{\alpha_{\sigma_n}} (t).
   $$
   Wir definieren den $\frak{k}$-Index 
   $   \frak{k} (\pmb{\alpha}) =   \frak{k} (\alpha_1,...,\alpha_n)$
   durch die folgende Gleichung
$$
   \frak{k} (\pmb{\alpha}) =\max \{k:
   \text{es existieren verschiedene Permutationen}\,\, \pmb{\sigma}_1,...,\pmb{\sigma}_k
$$
$$\hskip+2.5cm
\text{ mit der Eigenschaft:}\,\,\,
\forall \, j \,\,
 \forall t_0>0\,\, \exists\,\, t>t_0 \,\,
\text{mit}\,\,  \pmb{\sigma} (t) = \pmb{\sigma}_j \}.
$$

Insbesondere besagt
der Satz von Kan und Moshchevitin  \cite{KM}, dass
  die Gleichung
$$
 \frak{k} (\alpha_1,\alpha_2) =2
$$
für  inkommensurable Zahlen $\alpha_1,\alpha_2$
immer gilt. Es ist klar, dass
wir 
 für  fast alle $\pmb{\alpha}$
(im Sinne des Lebesguemaßes in $\mathbb{R}^n$)    $   \frak{k} (\pmb{\alpha}) = n!$ haben.

\vskip+0.3cm
Nun formulieren wir die Hauptergebnisse dieses Artikels.

\vskip+0.3cm

{\bf Satz 1.}\, {\it
Sei  $\pmb{\alpha} = (\alpha_1,...,\alpha_n)$  ein $n$-Tupel  paarweise inkommensurabler Zahlen.
   Dann  ist
   $$   \frak{k} (\pmb{\alpha}) \ge \sqrt{\frac{n}{2}}.$$}
   
 \vskip+0.3cm
 Für vier Zahlen gilt ein stärkeres, optimales Ergebnis  (siehe Satz 3 unten).
 
 Die folgende Behauptung
 zeigt,
 dass die Ungleichung von Satz 1 die optimale Größenordnung angibt.
  \vskip+0.3cm
 
{\bf Satz 2.}\, {\it
Sei $k\ge 3 $ und $ n = \frac{k(k+1)}{2}$.
Dann gibt es ein
  $n$-Tupel  $\pmb{\alpha}$  paarweise inkommensurabler Zahlen mit 
  $$   \frak{k} (\pmb{\alpha}) = k.$$}
  
 \vskip+0.3cm
 
 {\bf 2. Über Kettenbrüche.}
  \vskip+0.3cm
  
    Wir stellen jede irrationale Zahl $\alpha$ als  Kettenbruch 
  $$
  \alpha =[a_0;a_1,a_2,...,a_\nu,...],\,\,\,\,\, a_0 \in \mathbb{Z};\,\,\,a_j \in \mathbb{Z}_+,\, j =1,2,3,...  
  $$
  dar.
 Seien
  $$
  \frac{p_\nu}{q_\nu} = [a_0;a_1,a_2,...,a_\nu]
  $$
   die Näherungsbrüche.
  Wir verwenden die Notation
  $$
  \alpha_\nu = [a_\nu;a_{\nu+1},a_{\nu+2},...],\,\,\, 
 \alpha_\nu^* = [0;a_\nu,a_{\nu-1},a_{\nu-2},..., a_1] 
 \,\,\,\,
  \text{und}\,\,\,\,
  \xi_\nu = |q_\nu\alpha -p_\nu|.
 $$
 Es ist klar, dass die Formeln
\begin{equation}\label{inv}
\alpha_\nu^*=
 \frac{q_{\nu-1}}{q_{\nu}} , 
\end{equation}
\begin{equation}\label{zero}
\xi_\nu =
\frac{1}{q_\nu(\alpha_{\nu+1}+\alpha_\nu^*)}=
 \frac{1}{q_{\nu+1}+\frac{q_\nu}{\alpha_{\nu+2}}}
\end{equation}
und
$$
\psi_\alpha(t) = \xi_\nu,\,\,\,\,\, q_\nu\le t<q_{\nu+1}
$$
gelten.

\vskip+0.3cm

{\bf Hilfssatz 1.}\, {\it
Seien $ \alpha, \beta \in \mathbb{R}$ und $\frac{p_\nu}{q_\nu}, \frac{s_\mu}{r_\mu}$  Näherungsbrüche für $\alpha$ und $\beta$.
Dann  ist $ \alpha_{\nu}^* =   \beta_{\mu}^*$, 
falls
$
q_{\nu-1} =r_{\mu-1}, \,\,\, q_\nu =r_{\mu}
$
 ist.
}

\vskip+0.3cm
{\bf 
Beweis}. Die Behauptung folgt aus den Formeln (\ref{inv}) f\"{u}r $\alpha_{\nu} = \frac{q_{\nu-1}}{q_\nu}, \beta_{\mu} =  \frac{r_{\mu-1}}{r_\mu}$ und aus der Eindeutigkeit der Darstellung der
rationalen Zahlen als Kettenbr\"{u}che.$\Box$

\vskip+0.3cm
 
 {\bf 3. Der  Haupthilfssatz.}
  \vskip+0.3cm
  Für $\alpha, \beta \in \mathbb{R}$ verwenden wir die Notation von Hilfssatz 1,
  $ \xi_\nu = |q_\nu\alpha - p_\nu|, \eta_\mu = |r_\mu\beta - s_\mu|$.
  \vskip+0.3cm

{\bf Haupthilfssatz.}
\,
{\it
Sei $ d\ge 1$.
Nehmen wir an, dass die Ungleichungen
\begin{equation}
\label{1}
\xi_\nu\le \eta_\mu,
\end{equation}
\begin{equation}
\label{2}
\xi_{\nu+1}\le \eta_{\mu+d-1},
\end{equation}
\begin{equation}
\label{3}
q_{\nu+1} \le   r_{\mu+1},
\end{equation}
\begin{equation}
\label{4}
q_{\nu+2} = r_{\mu+d}
\end{equation}
gelten.
Dann werden die  Ungleichungen in (\ref{1}), (\ref{2}) und (\ref{3})   zu Gleichungen,
es gilt
\begin{equation}\label{und}
d = 2\,\,\,\,\text{und}
\,\,\,\,
\alpha_{\nu+2}^* = \beta_{\mu+2}^*.
\end{equation}
}

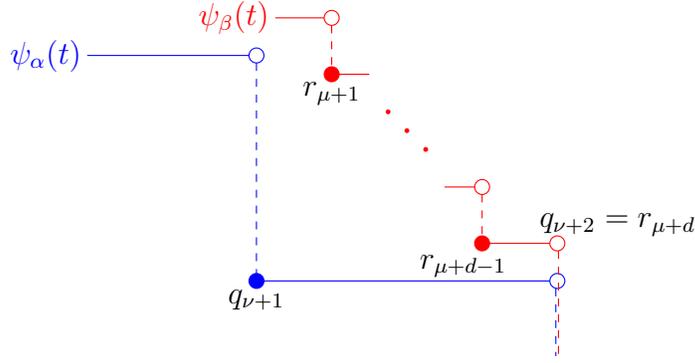
\begin{figure}[h]
  \centering
  \begin{tikzpicture}[scale=0.5]

 \draw[color=blue] (0,0) -- (7.8,0);

    \node[draw=blue,fill=blue,circle,inner sep=2pt] at (0,0) {};
    \node[draw=blue,circle,inner sep=2pt] at (8,0) {};

      \node[draw=blue,circle,inner sep=2pt] at (0,6) {};

 \draw (9.6,1.5)  [color=black] node {$q_{\nu+2} = r_{\mu+d}$};
 
 \draw (5.5, 0.4)  [color=black] node {$r_{\mu+d-1}$};

  \draw (2, 5)  [color=black] node {$r_{\mu+1}$};
  
   \draw (0, -0.5)  [color=black] node {$q_{\nu+1}$};

       \draw[color=red] (5,2.5) -- (5.8,2.5); 
        \node[draw=red,circle,inner sep=2pt] at (6,2.5) {};
     \draw[dashed,color=red]  (6,2.3) -- (6,1.1);
    
        \node[draw=red,fill=red,circle,inner sep=2pt] at (6,1) {};
         
      \draw[color=red] (6,1) -- (7.8,1);

      \draw[dashed,color=blue]  (7.96,-0.2) -- (7.96,-2);
      \draw[dashed,color=red]  (8.03,0.8) -- (8.03,-2);
        
    \node[draw=red,circle,inner sep=2pt] at (8,1) {};

      \draw[dashed,color=blue]  (0,5.8) -- (0,0);

 \draw (-5.6,6 )  [color=blue] node {$\psi_\alpha(t)$};

   \draw[color=blue] (-4.5,6) -- (-0.2,6);
  
  \draw[color=red] (0.5,7) -- (1.8,7);
  \node[draw=red,circle,inner sep=2pt] at (2,7) {};
        \draw[dashed,color=red]  (2,6.8) -- (2,5.5);
           \node[draw=red,fill=red,circle,inner sep=2pt] at (2,5.5) {};

  \node[draw=red,fill=red,circle,inner sep=0.5pt] at (3.5,4.5) {}; 
    \node[draw=red,fill=red,circle,inner sep=0.5pt] at (4,4) {}; 
      \node[draw=red,fill=red,circle,inner sep=0.5pt] at (4.5,3.5) {};

        \draw[color=red] (2,5.5) -- (3,5.5); 
    
     \draw (-0.6,7 )  [color=red] node {$\psi_\beta(t)$};

    \end{tikzpicture}
  \label{fig1}
    \caption{zum Haupthilfssatz: folgende Situation tritt nicht ein}
  \end{figure}
  
  \vskip+0.3cm

{\bf 
Beweis}.
F\"{u}r  $d=1$  ist dies  Hilfssatz 2 aus \cite{KM}. Tats\"{a}chlich zeigt
dieser Hilfssatz, dass   (\ref{1}) und (\ref{4})  gleichzeitig unmöglich sind.

Sei $ d \ge 2$ und 
 (\ref{1}) - (\ref{4})  mögen gelten.

1) Wir beweisen, dass
\begin{equation}
\label{a1}
\beta_{\mu+d+1} \ge \alpha_{\nu+3}.
\end{equation}
Aus (\ref{2}) und (\ref{zero}) folgt
$$
 \frac{1}{r_{\mu+d}+\frac{r_{\mu+d-1}}{\beta_{\mu+d+1}}} = \eta_{\mu+d-1} \ge \xi_{\nu+1} =
 {
  \frac{1}{q_{\nu+2}+\frac{q_{\nu+1}}{\alpha_{u+3}}}
  }
\,\,\,\,\, \text{
oder nach (\ref{4}),}
\,\,\,\,\,\,
\frac{q_{\nu+1}}{\alpha_{\nu+3}}\ge \frac{r_{\mu+d-1}}{\beta_{\mu+d+1}}.
$$
Somit ist
$$
\beta_{\mu+d+1} \ge \frac{r_{\mu+d-1}}{q_{\nu+1}} \, \alpha_{\nu+3} \ge \alpha_{\nu+3}
$$
nach (\ref{3}). Die Ungleichung 
(\ref{a1}) ist bewiesen.

Wir bemerken, dass 
wir 
statt (\ref{a1})   eine strikte Ungleichung
\begin{equation}
\label{a1f}
\beta_{\mu+d+1} > \alpha_{\nu+3}
\end{equation}
haben,
falls  die strikte Ungleichung 
\begin{equation}\label{mooree}
q_{\nu+1} <r_{\mu+d-1}
\end{equation}
gilt.

2) Wir betrachten  die Kontinuante
$$
b = \langle b_{\mu+2},b_{\mu+3},..., b_{\mu+d}\rangle,
 \,\,\,\,
 b_-= \langle b_{\mu+3},...,  b_{\mu+d}\rangle,\,\,\,
 \text{und}\,\,\,\,
b_-^-= \langle b_{\mu+3},...,  b_{\mu+d-1}\rangle.
$$
Dann ist
\begin{equation}\label{kont}
 r_{\mu+d}=  r_{\mu+1}b +   r_\mu b_- >  r_{\mu+1}b
\end{equation}
und nach dem Perron'schen Satz hat man
\begin{equation}\label{pe}
 \beta_{\mu+2} b_- = b + \frac{(-1)^d}{ b_- \beta_{\mu+d+1}+{b_-^-}}\le b + \frac{1}{ \beta_{\mu+d+1}}.
\end{equation}

Wir beweisen, dass
\begin{equation}
\label{a}
  b \le a_{\nu+2}.
\end{equation}

Falls $a_{\nu+2} \le b-1$ ist,  haben wir
die Ungleichung
$$
q_{\nu+2} = a_{\nu+2}q_{\nu+1} + q_\nu \le (b-1) q_{\nu+1} + q_\nu< bq_{\nu+1}\le b r_{\mu+1} < r_{\mu+d}
$$
aus (\ref{3}) und (\ref{kont}).  
Das liefert einen Widerspruch zwischen der letzten Ungleichung und (\ref{4}). Die Ungleichung (\ref{a}) ist bewiesen.

3) 
Wir bemerken, dass die Ungleichung   
\begin{equation}
\label{b}
\beta_{{\mu}+2} b_-  \le\alpha_{\nu+2}
\end{equation}
aus (\ref{pe},\ref{a})  und (\ref{a1})
folgt.
Tatsächlich ist
$$
\beta_{{\mu}+2} b_- \le
b + \frac{1}{ \beta_{\mu+d+1}}\le a_{\nu+2}+ \frac{1}{\alpha_{\nu+3}} = \alpha_{\nu+2}
$$
und (\ref{b}) ist bewiesen. Falls (\ref{mooree}) gilt,  folgt    die strikte Ungleichung
\begin{equation}
\label{bf}
\beta_{{\color{red}\mu}+2} b_-  <\alpha_{\nu+2}
\end{equation}
aus (\ref{a1f}).

4) 
{
Betrachten wir den Wert
$$
x = r_{\mu+1} - q_{\nu+1}.
$$
Wir
beweisen, dass die Ungleichung}
\begin{equation}\label{w}
0\le bx \le q_{\nu}-r_\mu b_-
\end{equation}
gilt.
Betrachten wir den Wert
$$
x = r_{\mu+1} - q_{\nu+1}.
$$
Nach (\ref{3}) hat man $ x\ge 0$.
 Aus (\ref{1}) und (\ref{zero}) folgt
$$
\frac{1}{r_{\mu+1}+\frac{r_\mu}{\beta_{\mu+2}}} = \eta_{\mu} \ge \xi_{\nu} = \frac{1}{q_{\nu+1}+\frac{q_{\nu}}{\alpha_{\nu+2}}},
$$
oder
$$
x    \le  \frac{q_{\nu}}{\alpha_{\nu+2}} - \frac{r_\mu}{\beta_{\mu+2}} \le \frac{q_{\nu}-r_\mu b_-}{\alpha_{\nu+2}} \
,
$$
wegen (\ref{b}). 
Aus  (\ref{a})
erhält  man
$$
bx \le \alpha_{\nu+2} x \le q_{\nu}-r_\mu b_-,
$$
 und (\ref{w}) folgt daraus.
 
 Falls (\ref{mooree}) gilt,   haben wir aus (\ref{bf}) die Ungleichung
 \begin{equation}\label{zerr}
bx< q_{\nu}-r_\mu b_-.
 \end{equation}
  
 5) Wir bemerken, dass   die Ungleichung 
 $  q_{\nu+1}\le r_{\mu+1} \le r_{\mu+d-1} $
 aus (\ref{3})  folgt.
 Nun betrachten wir  die folgenden zwei Fälle.
 
 \vskip+0.3cm
 {\bf Fall 1$^0$.}  $  q_{\nu+1}= r_{\mu+d-1}$. Dann ist $r_{\mu+d-1}\le r_{\mu+1}$ und $ d =2$, und es gilt
 $$
 q_{\nu+1} = r_{\mu+1},\,\,\,\,\, q_{\nu+2} = r_{\mu+2}.
 $$
 Sodass (\ref{und})  nach 
Hilfssatz 1 gilt .
 
 \vskip+0.3cm
  {\bf Fall 2$^0$.} Die Ungleichung (\ref{mooree}) gilt.
  Daraus { folgt die} strikte Ungleichung  (\ref{zerr}).
   Sodass
\begin{equation}\label{www}
0\le bx < q_{\nu}-r_\mu b_-
\end{equation}
 statt (\ref{w}) gilt.
\vskip+0.3cm
 
  Wir bemerken, dass
 $$
 q_{\nu+2} = a_{\nu+2} q_{\nu+1} + q_\nu \equiv q_\nu \pmod{
 {
 q_{\nu+1}}
 }.
 $$
 Dann
 ist 
$$
bx = b(r_{\mu+1}- q_{\nu+1} )\equiv br_{\mu+1} = r_{\mu+d} - r_\mu b_- =
q_{\nu+2}- r_\mu b_- \equiv q_\nu - r_\mu b_- \pmod{q_{\nu+1}}
$$
wegen (\ref{kont}) und (\ref{4})
oder
\begin{equation}\label{q}
bx \equiv q_\nu - r_\mu b_- \pmod{q_{\nu+1}}.
\end{equation}
Da $ 0< q_\nu - r_\mu b_- < q_\nu <q_{\nu+1}$ ist,
 erhalten wir einen Widerspruch zwischen  (\ref{www}) und (\ref{q}).

Der Hilfssatz ist  somit bewiesen.$\Box$

\vskip+0.3cm

{\bf 4. Beweis von Satz 1.}

\vskip+0.3cm
Für $\alpha_i$ betrachten wir  den Kettenbruch
$\alpha_i  = [a_{i,0}; a_{i,1},...,a_{i,\nu},...]$.
Bezeichnen wir
$\alpha_{i,\nu}^* = [0; a_{1,\nu},...,a_{i,1}]
$.

Sei $\frak{T}$ die Menge der Punkte, an denen
   mindestens eine der Funktionen 
$\psi_{\alpha_j }(t), 1\le j \le n$
nicht stetig ist. 
Für $T\in \frak{T}$ betrachten wir den Wert
$$
\tau (T) = |\{ j \in \{1,...,n\}:\,\, \psi_{\alpha_j}(t)\,\,
\text{ist nicht stetig im Punkt } \, T\}|
$$

Sei
$$ W = \frac{\sqrt{8n+1}-1}{2}.$$
Es gibt    zwei Fälle.

\vskip+0.3cm
{\bf Fall 1$^0$.} für alle hinreichend großen $T$ aus $\frak{T}$ gilt  $ \tau (T) \le W$;

{\bf Fall 2$^0$.} 
es gibt unendlich viele 
 $T\in \frak{T}$ mit $ \tau (T) > W$.
 \vskip+0.3cm
 
 Wir betrachten den  ersten Fall. Sei $T$ groß genug. Nehmen wir an, dass
$\pmb{\sigma} (T) = (1,2,3,...,n)$ ist, sodass 
 $$
 \psi_{\alpha_{1}} ({T}) >
 \psi_{\alpha_{2}} ({T}) >
  \psi_{\alpha_{3}} ({T}) >
  ... >
   \psi_{\alpha_{n}} ({T}).
   $$
 Wir definieren die Werte  
$$
T_0=T,\,\,\,\,\,
T_j = \min\{t\ge T : \psi_{\alpha_j}(t)<  \psi_{\alpha_{j+1}}(t)\}, \,\,\, 1\le j \le n-1,
$$
die nach dem Satz von Kan und Moshchevitin  existieren.
Dann ist  die Funktion $\psi_{\alpha_j} (t)$ im Punkt $T_j$ nicht stetig
und es gelten die Ungleichungen
 \begin{equation}\label{fur0}
 \psi_{\alpha_j} (T_j)< \psi_{\alpha_{j+1}} (T_j) 
 \end{equation}
und
 \begin{equation}\label{fur}
 \psi_{\alpha_j} (t)> \psi_{\alpha_{j+1}} (t)\,\,\,\,
 \text{ für alle}\,\, t \in [T,T_j).
 \end{equation}
 Aus (\ref{fur0}) und (\ref{fur}) sehen wir, dass die Permutation 
 $\pmb{\sigma}(T_j)$  sich von anderen Permutationen 
 $\pmb{\sigma}(t)$ mit $ t \in [T_0, T_j)$   unterscheidet. 
 
Nun ordnen wir die Punkte  $T_j, 0\le j \le n-1 $:
 \begin{equation}\label{fur0o}
 T_0\le T_{j_1}\le T_{j_2}\le T_{j_3}\le ...\le T_{j_{n-1}}.
 \end{equation}
 Sei $ z$ die Anzahl der strikten  Ungleichungen in (\ref{fur0o}).
 Aus der Annahme des Falles   haben wir
 $(z+1) W\ge n$. Nun haben wir   mindestens $ (z+1)\ge n/W$ verschiedene Werte $T_j$ und 
 damit  $ \ge n/W$ verschiedene Permutationen $\pmb{\sigma}(T_j)$.
  Somit folgt
  \begin{equation}\label{02}
    \frak{k} (\pmb{\alpha})\ge \frac{n}{W} \ge \sqrt{\frac{n}{2}}.
  \end{equation}
 
 \vskip+0.3cm
 Nun 
  betrachten wir den  zweiten Fall.  Wir nehmen   ein ausreichend  großes  $T \in \frak{T}$ mit $ \tau(T) > W$.
  Wir betrachten die Zahlen $\alpha_j$, für die die Funktionen 
  $\psi_{\alpha_j}(t)$ im Punkt $T$ nicht  stetig sind.
  Es gibt $\ge [W]+1$  dieser Zahlen. 
Wir bezeichnen   diese  Zahlen  als $ \alpha_1,..., \alpha_{W_1},\, W_1\ge [W]+1$ und
betrachten   die Werte 
$$
\psi_{\alpha_j} (T-1) =   \lim_{t\to T-} \psi_{\alpha_j} (t).
$$
Ohne Beschränkung der Allgemeinheit nehmen wir an, dass
  $$
   \psi_{\alpha_1} (T-1) >   \psi_{\alpha_2} (T-1)> \psi_{\alpha_3} (T-1)>...>   \psi_{\alpha_{W_1}} (T-1)
$$
ist.
  Es ist klar, dass
  für jedes $ j\in \{2,...,W_1\}$  und für alle $t\in [T-1,T)$ 
  die Ungleichung
 $$
 \psi_{\alpha_j} (t) < \min_{i\in \{1,...,j-1\}} \psi_{\alpha_i}(T )= \min_{i\in \{1,...,j-1\}} \psi_{\alpha_i}(t) 
 $$
 gilt.
  Sei $R_1 = T$.
  Wir betrachten die Werte
 \begin{equation}\label{li0}
  R_j= {\sup \{ t<T:  \,\, 
   \psi_{\alpha_j} (t) > \min_{i\in \{1,..., j-1\}} \psi_{\alpha_i}(t)\} <T} ,\,\,\,\, j \ge 2.
\end{equation}
Es folgt  nach dem Haupthilfssatz, dass
  \begin{equation}\label{li}
  \psi_{\alpha_j} (R_j-1) =
    \lim_{t\to R_j} \psi_{\alpha_j} (t)>   \max_{i \in \{1,...,j-1\}} \,\, \psi_{\alpha_j} (R_j)
  \end{equation}
  gilt.
  Falls  { nämlich}
  $$
  \lim_{t\to R_j} \psi_{\alpha_j} (t)\le  \psi_{\alpha_k} (R_j),\,\,\, \text{mit}\,\,\,  1\le k \le j-1
  $$  
  ist, so
  betrachten wir die Indizes $ \mu,\nu$ mit 
  $q_\nu = r_\mu = T$ (hier ist $q_\nu$ ein Teilnenner des Kettenbruches für $\alpha_j$ und 
    $r_\nu$ ein Teilnenner des Kettenbruches für $\alpha_k$). Nun folgt  aus dem Haupthilfssatz, dass
   $ \alpha_{j,\nu}^* = \alpha_{k,\mu}^*$ gilt, was    für großes $T$ unmöglich ist.
   
   Aus (\ref{li}) und der Definition (\ref{li0}) des Wertes $R_j$ sehen wir, dass
   die Permutation  $\pmb{\sigma} (R_j-1)$ sich
     von allen anderen Permutationen  $\pmb{\sigma} (t)$ mit $ t \in [R_j, T)$   unterscheidet. So ist
 \begin{equation}\label{bei}
   \frak{k} (\pmb{\alpha})
  \ge\,\,
  \text{die Anzahl  verschiedener Zahlen zwischen}\,\,\,
  R_1, R_2,...,R_{W_1}.
  \end{equation}

Nun
beweisen wir, dass
$$
R_i = R_j = R_k \,\, \text{mit verschiedenen}\,\, i,j,k \in \{2,...,W_1\}
$$
unmöglich ist.
Seien
$$ 
R_i = q_{\nu},\,\,\, R_j = r_\mu ,\,\,\, R_k = t_\lambda
$$
die Teilnenner der Kettenbrüche für $\alpha_j, \alpha_i, \alpha_k$, bzw.
Wir nehmen
ohne Beschränkung der Allgemeinheit  an, dass
$$
\psi_{\alpha_i} (R_i) < \psi_{\alpha_j} (R_j)< \psi_{\alpha_k} (R_k)
  $$
  gilt.
  Falls $ T = q_{\nu+1} = r_{\mu+1}$ ist, folgt dann $ \alpha_{k,\nu+1}^* = \alpha_{i,\mu+1}^*$ nach 
   Hilfssatz 1,
  was    für großes $T$  unmöglich ist.
  So ist entweder $ T  = q_{\nu+d}$ oder $ T = r_{\mu+d}$ mit $ d \ge 2$. 
  Sei $ T = t_{\lambda+ d'}$ ein Teilnenner des Kettenbruches für $\alpha_k$.
  Dann   haben wir
  entweder $ \alpha_{j,\nu+d}^* = \alpha_{k,\lambda+d'}^*$ oder 
  $ \alpha_{i,\mu+d}^* = \alpha_{k,\lambda+d'}^*$ nach dem Haupthilfsatz. Diese beiden Fälle sind  unmöglich.
  
  Somit gibt es zwischen $R_1, R_2,...,R_{W_1}$ 
  mindestens $\left[ \frac{W_1}{2}\right]$  verschiedene Zahlen und nach  (\ref{bei}) gilt
$$
    \frak{k} (\pmb{\alpha})\ge
  \left[ \frac{W_1}{2}\right]\ge \left[ \frac{[W]+1}{2}\right]\ge \frac{W+1}{2} =\frac{n}{W} \ge \sqrt{\frac{n}{2}}.
  $$

Daraus folgt die Behauptung.

\vskip+0.3cm

{\bf 5. Vier Zahlen.}
\vskip+0.3cm
Hier beweisen wir eine stärkere Behauptung im Fall $n=4$. Aus  Satz 2
folgt, dass die Ungleichung von Satz  3 optimal ist.

\vskip+0.3cm

{\bf Satz 3.}\, {\it Seien  $ \alpha_i ,  1\le i \le 4$ vier irrationale Zahlen mit
\begin{equation}\label{ijn}
\alpha_i \pm \alpha_j \not\in \mathbb{Z},\,\,\,\,\, i \neq j.
\end{equation}
Dann ist
$$
\frak{k} (\alpha_1,\alpha_2,\alpha_3,\alpha_4) \ge 3.
$$}

\vskip+0.3cm

{\bf 
Beweis}. \,
Falls  $\frak{k} (\alpha_1,\alpha_2,\alpha_3,\alpha_4) = 2$,
können wir annehmen, dass
für hinreichend großes
$t$ entweder
\begin{equation}\label{s1}
\psi_{\alpha_1}(t)\ge 
\psi_{\alpha_2}(t)\ge 
\psi_{\alpha_3}(t)\ge 
\psi_{\alpha_4}(t)
\end{equation}
oder
\begin{equation}\label{s2}
\psi_{\alpha_4}(t)\ge 
\psi_{\alpha_3}(t)\ge 
\psi_{\alpha_2}(t)\ge 
\psi_{\alpha_1}(t)
\end{equation}
gelten. Wir nehmen an, dass (\ref{s1}) f\"{u}r $q_0\le t < q_1$ gelten und 
(\ref{s2}) f\"{u}r $q_1\le t < q_2$ gelten.
Dann   sind die Funktionen 
$
\psi_{\alpha_1}(t), 
\psi_{\alpha_2}(t), 
\psi_{\alpha_3}(t)
$ im Punkt $q_0$ 
nicht stetig,
und 
  die Funktionen 
$
\psi_{\alpha_2}(t), 
\psi_{\alpha_3}(t), 
\psi_{\alpha_4}(t)
$ im Punkt $q_1$ 
nicht stetig. Daher   ist $q_1$ ein gemeinsamer Nenner der Näherungsbrüche für $\alpha_1, \alpha_2,\alpha_3$,
und 
  $q_2$ ein gemeinsamer Nenner der Näherungsbrüche für $\alpha_2, \alpha_3,\alpha_4$.
Wir bezeichnen die  Nenner  der Näherungsbrüche  für $\alpha_i$ mit
  $ q^{(i)}_{l}, l =0,1,2,3,... $\, . 
Es sind
$$
q_{1} = q_{\lambda}^{(1)}  = q_{\nu}^{(2)} = q_{\mu}^{({\color{red}3})},\,\,\,\,\,
q_{2} = q_{{\color{red}\nu}'}^{(2)}  = q_{{\color{red}\mu}'}^{(3)} = q_{\kappa'}^{(4)}
$$
mit $ \lambda, \nu,\nu'\mu,\mu',\kappa' \in \mathbb{Z}_+$.
Es ist klar, dass $ q_1, q_2$   gemeinsame Nenner f\"{u}r $\alpha_2$ und $\alpha_3$ sind  und
$$
\psi_{\alpha_2} (t) \le  \psi_{\alpha_3}(t),\,\,\,\,\ q_1 \le t < q_2
$$
gilt.

Betrachten wir die folgende Fälle.

{\bf 1$^0$.}
Falls $ \nu' = \nu+1, \mu' = \mu+1$ ist,
folgt die Gleichung $ \alpha_{2,\nu+1}^* = \alpha_{3,\mu+1}^*$
aus   Hilfssatz 1.

{\bf 2$^0$.}
Im  Fall $  \nu' > \nu+1, \mu' = \mu+1$ 
haben wir die Gleichung
$ \alpha_{2,\nu'}^* = \alpha^*_{3,\mu'}$
nach dem  Haupthilfssatz.

{\bf 3$^0$.} Falls $ \mu' > \mu+1$ gilt, betrachten wir die Funktionen
$$
\psi_{\alpha_3} (t) \le  \psi_{\alpha_4}(t),\,\,\,\,\ q_1 \le t < q_2.
$$
Wir bemerken, dass $q_2 = q_{\mu'}^{(3)} = q_{\kappa'}^{(4)}$ gilt.
 Dann 
 folgt
  die Gleichung
$ \alpha_{4,\kappa'}^* = \alpha_{3,\mu'}^*$
 nach dem Haupthilfssatz.

Nun sehen wir, dass  es 
 die Indizes  $i,j,\, i\neq j$   und unendlich viele $\mu,\nu$ mit
$$
\alpha_{i,\mu}^* = \alpha_{j,\nu}^*
$$
gibt.
Das bedeutet,  dass $ \alpha_i \pm \alpha_j \in \mathbb{Z}$ gilt
und
alles ist bewiesen.$\Box$

\vskip+0.3cm

 {\bf Bemerkung.}
Hier sei angemerkt, dass Sergei Konyagin \cite{ko} uns ein einfaches, elegantes Argument
mitgeteilt hat, wie man aus dem Satz  von  Erd\"{o}s  und Szekeres \cite{ES}
die folgende Ungleichung  aus  Satz 3 bekommen kann:
{\it falls für
$\pmb{\alpha} = {(\alpha_1 ,..., \alpha_n)} \in (\mathbb{R}\setminus\mathbb{Q})^n$ 
die Bedingung (\ref{ijn}) gilt,
 folgt
$ \frak{k} (\pmb{\alpha}) \ge \gamma \log\log n$ mit einer absoluten positiven Konstante $\gamma$}.

\vskip+0.3cm

{\bf 6. Beweis von Satz 2.}

   \vskip+0.3cm

  Zuerst nummerieren wir
 die Menge der 
 $ n=  \frac{k(k+1)}{2}
$ Elemente des $n$-Tuples
als
$$
\frak{A}_k = \{ \frak{a}_1,\frak{a}_2,...,\frak{a}_k, \frak{a}_{k+1}, \frak{a}_{k+2},...,\frak{a}_{n-1},\frak{a}_n\}=
$$
$$
=
\left\{ 
\{1\}, \underbrace{\{1,k\}, \{1,k-1\},....,\{1,2\}}_{k-1\,\, \text{Elemente}};
\{2\}, \underbrace{\{2,k\},....,\{2,3\}}_{k-2 \,\,\text{Elemente}}; \{3\}, 
\cdots ; \{k-1\} ,\underbrace{\{k-1,k\}}_{1 \,\,\text{Element}}; \{k\}
\right\}.
$$
Für
$$
w_1 = 1,\,\,\,\, w_{i+1} = w_i+k-i+1,\,\,\,\, i= 1,...,k
$$
haben wir
$$
\frak{a}_{w_i}  =  \{i\},\,\,\, 1\le i \le k;\,\,\,\,\,
 \frak{a}_{w_i +j }  =\{ i,k-j+1\}
 ,\,\, 1\le j \le k-i,
  \,\,\,\text{oder}\,\,\,
  \frak{a}_{w_i +k-j+1 }  =\{ i,j\},\,\, i+1\le j \le k
$$
und
für
$$
\frak{X} =
\left\{ 
\{x_1\}, \underbrace{\{x_{1,2}, y_{1,2}\}, \{x_{1,3}, y_{1,3}\},....,\{x_{1,k}, y_{1,k}\}}_{k-1\,\, \text{Elemente}};
\{x_2\},  
\cdots ;
\{x_{k-1}\}, \underbrace{\{x_{k-1,k}, y_{k-1,k}\}}_{1 \,\,\text{Element}}; \{x_k\}
\right\}
$$
definieren wir
$$
\underline{\frak{X}} = \underline{\frak{X}}_1 =
\left\{ 
\{x_1\}, \underbrace{\{x_{1,2}, y_{1,2}\}, \{x_{1,3}, y_{1,3}\},....,\{x_{1,k}, y_{1,k}\}}_{k-1\,\, \text{Elemente}}
\right\}
$$
und
$$
  \underline{\frak{X}}_i =
\left\{ 
\{x_i\}, \underbrace{\{x_{i,i+1}, y_{i,i+1}\}, \{x_{i,i+2}, y_{i,i+2}\},....,\{x_{i,k}, y_{i,k}\}}_{k-i\,\, \text{Elemente}}
\right\},\,\,\, i=2,..., k,
$$
sodass
$
\frak{X} =\left\{
\underline{\frak{X}}_1; \underline{\frak{X}}_2;\cdots; \underline{\frak{X}}_k
\right\}
$ ist.
Dann ist
$$
  \underline{\frak{A}_k}_i =
  \left\{ \{i\},
 \underbrace{   \{i,k\},\{i,k-1 \} , \{i,k-2 \},....,\{i,i+2 \}, \{i,i+1 \}}_{k-i\,\, \text{Elemente}}
\right\},\,\,\, i=1,..., k.
  $$

Wir betrachten die Permutation $\Sigma$, die durch die Gleichungen
\begin{equation}\label{gaze0}
\begin{gathered}
\Sigma(\{1\}) = \{k\},
\\
\Sigma(\{i\}) = \{i-1\},\,\,\, 2\le i \le k,
\\
\Sigma(\{ 1,j \}) = \{ j-1,k\},\,\,\, 2\le j \le k,
\\
\Sigma(\{i,j \}) = \{i-1,j-1\},\,\,\, 2\le i \le k,\,\,\, i+1\le j \le k 
\end{gathered}
\end{equation}
definiert ist.

\vskip+0.3cm
{\bf Hilfssatz 2.} {\it
 
{\rm (i)}
 Für jedes $ l =  0,..., k-1$ ist
 $$
 \underline{\Sigma^{l}(\frak{A}_k)} =
\left\{ \{l+1\}, 
 \underbrace{\{l, l+1\}, \{l-1, l+1\},....,\{1, l+1\}}_{l\,\, \text{Elemente}},
 \underbrace{\{l+1,k\}, \{l+1, k-1\},....,\{l+1, l+2\}}_{k-l-1\,\, \text{Elemente}}\right\}.
$$

{\rm (ii)}
Die Ordnung  der Permutation $\Sigma$ ist $k$.
}

\vskip+0.3cm
{\bf 
Beweis.
}
Wir
betrachten  $\frak{X} = \Sigma^{l}(\frak{A}_k)$ und
$
\underline{\Sigma^{l}}_i = \underline{\frak{X}}_i. 
$
Dann folgt
$$
\underline{\Sigma^{l}}_{i}  =
 \left\{ \{l+i\}, 
 \underbrace{\{l, l+i\}, \{l-1, l+i\},....,\{1, l+i\}}_{l\,\, \text{Elemente}},
 \underbrace{\{l+i,k\}, \{l+i, k-1\},....,\{l+i, l+i+1\}}_{k-l-i\,\, \text{Elemente}}\right\}
$$
$$
\text{
für}\,\,\, 1\le i \le k-l
$$ 
und
$$
\underline{\Sigma^{l}}_{i}  =
 \left\{ \{ l+i-k\}, 
 \underbrace{\{l+i-k,l\}, \{l+i-k,l-1\},....,\{l+i-k, l+i-k+1\}}_{k-i\,\, \text{Elemente}}\right\}
$$
$$
\text{
für}\,\,\,  k-l+1\le i \le k
$$
aus (\ref{gaze0})
mittels Induktion nach $l\ge 0$. Beide Aussagen (i) und (ii) folgen daraus.$\Box$

\vskip+0.3cm

Zum Beispiel  gilt für  $k=3, n = 6$

\begin{equation}\label{gaze}
\begin{gathered}
 \frak{A}_3 =\Sigma^3(\frak{A}_3)  = \{ \{1\},\{1,3\},\{1,2\};\{2\}, \{2,3\};\{3\}\},\\
 \Sigma(\frak{A}_3) = \{ \{2\},\{1,2\},\{2,3\};\{3\}, \{1,3\};\{1\}\},\\
 \Sigma^2(\frak{A}_3) =\{ \{3\},\{2,3\},\{1,3\};\{1\}, \{1,2\};\{2\}\}
\end{gathered}
\end{equation}
und
$$
\underline{\frak{A}_3} =\{ \{1\},\{1,3\},\{1,2\}\},\,\,\,\,
 \underline{\Sigma(\frak{A}_3)}=  \{ \{2\},\{1,2\},\{2,3\}\},\,\,\,\,
  \underline{\Sigma^2(\frak{A}_3)}= \{ \{3\},\{2,3\},\{1,3\}\}
$$
(siehe das Diagramm in Figur 2).

\vskip+0.3cm

Daher sind
\begin{equation}\label{permuta}
\Sigma^l(A_k)
=
\{ \frak{a}_1^l,\frak{a}_2^l,...,\frak{a}_n^l\}
,\,\,\,\, 0\le l \le  k-1
\end{equation}
verschiedene Permutationen  und $ \Sigma^k(A_k) = A_k$.

Um  die Zahlen $\alpha^{(i)}, 1\le i\le k$ und $ \alpha^{(i,j)},   i<j\le k$
 zu konstruieren,   brauchen wir eine induktive Prozedur.
Wir
 definieren   
 $$ 
 t_{-i} = 0,\,\,\,\,  1\le i \le k;\,\,\,\,\,
 t_{i} = q_1^{(i+1,j)} =q_0^{(i+1,j)}= q_0^{(i+1)}= a_1^{(i+1,j)} = 1,\,\,
 0\le i  \le k-1,\ i+1<j\le k.
 $$
 Sei $\nu\ge 1$. 
Nehmen wir an, dass die Werte
$$
t_s,\,\,\,\,\,
s\le k\nu -1
$$
mit
\begin{equation}\label{schritt}
(t_{k(\nu-1)+i}, t_{k(\nu-2)+i})
=1,\,\,\, 0\le i\le k-1;\,\,\,\,\,
(t_{k(\nu-1)+i}, t_{k(\nu-1)+j}) =1,\,\,\,\, 0\le i {\color{red} <j}  \le k-1
\end{equation}
und die Teilnenner
$$
a^{(j)}_{s} ,\,\,  1\le s\le \nu-1,\, 1\le j \le k ;
 \,\,\,\,\,\,\,\,\,\,\,
a^{(i,j)}_{s} ,\,\,  1\le s\le 2\nu-1,\, 1\le i< j \le k
 $$
 schon
definiert sind.
Dann betrachten wir die simultanen Kongruenzen
  \begin{equation}\label{z0}
  \begin{cases}
   x\equiv t_{k(\nu-2)}\pmod{t_{k(\nu-1)}}
   \cr
    x\equiv t_{k(\nu-1)}\pmod{t_{k(\nu-1)+j}},\,\,  1\le j \le k-1
   \end{cases}
\end{equation}
und ihre  Lösung 
$$
t_{k\nu}=x, \,\,\, 
(t_{k\nu}, t_{k\nu-j})   = 1,\,\,\, 1\le j \le k
$$
mit  einem hinreichend großen $x$.
Seien
\begin{equation}\label{teil}
a^{(1)}_{\nu}= \frac{x-t_{k(\nu-2)}}{t_{k(\nu-1)}};\,\,\,\,\,\,\,\,\,\,\,
a^{(1,j+1)}_{2\nu}= \frac{x-t_{k(\nu-1)}}{t_{k(\nu-1)+j}},\,\,\,
 1\le j \le k-1.
\end{equation}
Damit haben wir $a_\nu^{(i)}$ und $a_{2\nu}^{(i,j)}$   fur alle
$ \{i\}, \{i,j\} \in \underline{\frak{A}_k}$ definiert.

   Dann betrachten wir nacheinander  für $ l=1,..., k-1$    
    die Kongruenzen
 \begin{equation}\label{z0l}
  \begin{cases}
   x\equiv t_{k(\nu-2)+l}\pmod{t_{k(\nu-1)+l}}
   \cr
    x\equiv t_{k(\nu-1)+l}\pmod{t_{k(\nu-1)+l + j}},\,\,  1\le j \le k-1
   \end{cases}.
\end{equation}

   Wir definieren   die Werte
   $
   t_{k\nu+l}=x
   $
   als  Lösungen von  (\ref{z0l}) mit    hinreichend großem $x$ . Nun definieren wir die Teilnenner
   \begin{equation}\label{teil1}
   a^{(l+1)} _{\nu} \!=\! \frac{x-t_{k(\nu-2)+l}}{t_{k(\nu-1)+l}};\,
   a^{(l+1,j+1)}_{2\nu}\!= \!\frac{x-t_{k(\nu-1)+l}}{t_{k(\nu-1)+j}},\, l\!+\!1\le j \le k\!-\!1;\,
   a^{(i,l+1)}_{2\nu+1}\!\!= \frac{x-t_{k(\nu-1)+l}}{t_{k(\nu-1)+k+i-1}},\,1\!\le i\! \le l
   ,
   \end{equation}
 die der Menge $ \underline{\Sigma^{l}(\frak{A}_k)}$ entsprechen.
  Wir sehen, dass die Kongruenzen  (\ref{z0l}) analog zu den  Kongruenzen (\ref{z0}) sind.
   Falls   wir $ l=0$ in 
    (\ref{z0l}) nehmen, dann erhalten wir (\ref{z0}). 
    Aber die Definitionen für  die Teilnenner $a^{(i,j)}_\nu$    in (\ref{teil1}) und in   (\ref{teil}) sind unterschiedlich.
    In (\ref{teil})  sind  alle  $a^{(1,j)}_{2\nu}$  die zweite Elemente der Periode der Konstruktion, hingegen
    in   (\ref{teil1})  
    sind
    die Teilnenner
    $a^{(l+1,j)}_{2\nu}$   
    die zweite Elemente der Periode, unterdessen sind
    $a^{(i,l+1)}_{2\nu+1}$
    die erste Elemente der Periode.

     In allen Fällen,  nehmen wir an, dass
    \begin{equation}\label{z00}
    \frac{t_{s+1}}{t_s}  \to \infty,\,\,\,\, s\to \infty
    \end{equation}
    (die Werten $x$ sind groß genug in aller Lösungen).
    Außerdem können wir  den Parameter $t_{k\nu+i}$ so wählen, dass
        \begin{equation}\label{dop}
    \frac{t_{k\nu+i}}{t_{k\nu+i-1}}\cdot \frac{t_{k(\nu-1)+i-1}}{t_{k(\nu-1)+i}}  \to \infty,\,\,\,\, \nu\to \infty,
    \end{equation}
    weil
    $ k\nu+i > \max( {k\nu+i-1}, {k(\nu-1)+i})$ gilt.
    
Hiermit sind 
die Zahlen
$$
t_{k\nu +i},\,\,\,\,  0\le i \le k-1
$$
und die Teilnenner
$$
a^{(j)}_{\nu} ,\,\, 1\le j \le k ;
 \,\,\,\,\,\,\,\,\,\,\,
a^{(i,j)}_{2\nu} ,\,\,   a^{(i,j)}_{2\nu+1} ,\,\, 1\le i< j \le k
 $$
definiert. Aus (\ref{schritt}) und (\ref{z0},\ref{z0l}) sehen wir, dass (\ref{schritt}) für den nächsten Schritt der induktiven Prozesse auch gilt.

Wir betrachten die $n$ irrationale Zahlen
$$
\alpha^{(s)}
=[0; a_1^{(s)},...,a_\nu^{(s)},...],\,\, 1\le s \le k;\,\,\,\,\,\,
\alpha^{(i,j)}
= 
[0; a_1^{(i,j)},...,a_\nu^{(i,j)},...],\,\, 1\le i<j \le k.
$$
Seien 
$$ 
\psi_s (t) = \psi_{\frak{a}_{w_s}} (t), \, 1\le s\le k;\,\,\,\,\,
\psi_{i,j} (t) = \psi_{ \frak{a}_{w_i +k-j+1 } }(t), \, 1\le   i<j\le k
$$
die zugehörigen
  Funktionen des Irrationalitätsmaßes.
Für die Nenner der Näherungsbrüche gelten die Gleichungen
\begin{equation}\label{TK}
q_\nu^{(s)} = t_{k\nu+s-1},\,\, 1\le s \le k;\,\,\,\,\,\,\,
q_{2\nu}^{(i,j)} =  t_{k\nu+i-1},
,\,\,  
q_{2\nu+1}^{(i,j)} =  t_{k\nu+j-1}
\,\, 1\le i<j\le k.
\end{equation}
Seien 
$$
\alpha^{(s)}_\nu
=[a_\nu^{(s)};a_{\nu+1}^{(s)}...],\,\, 1\le s \le k;\,\,\,\,\,\,
\alpha^{(i,j)}_\nu
= 
 [a_\nu^{(i,j)}; a_{\nu+1}^{(i,j)} ,...],\,\, 1\le i<j \le k.
$$

 Die Konstruktion von Satz 2  mit drei Permutationen  für sechs Zahlen ist in Figur 2 dargestellt.

\begin{figure}[h]
  \centering
  \begin{tikzpicture}[scale=0.7]
 
    \draw[color=red] (-1,9.1) -- (1,9.1);
    
    \draw[color=red] (-1.3,9.4)  node {$
 \psi_1$};
\draw[color=gray] (-1.1,8.4)  node {$
 \psi_{1,2}$};
    
     \draw[color=gray] (-0.6,8.5) -- (1,8.5);
     
     \draw[color=yellow] (-1.2,8.9)  node {$
 \psi_{13}$};
     \draw[color=yellow] (-0.75,8.8) -- (1,8.8);
     \draw[color=green] (-1.3,7.2)  node {$
 \psi_{2}$};
     
     \draw[color=green] (-0.9,7) -- (6,7);
     
     \draw[color=blue] (-1.1,6.6)  node {$
 \psi_{2,3}$};
     \draw[color=blue] (-0.6,6.6) -- (6,6.6);

    \draw[color=gray] (1,6.8) -- (6,6.8);

\draw[color=black] (-1.2,5.1)  node {$
 \psi_{3}$};

     \draw[color=black] (-0.8,5) -- (11,5);
     
     \draw[color=yellow] (1,4.6) -- (11,4.6);
     
      \draw[color=blue] (6,4.8) -- (11,4.8);
  
             \draw[color=red] (1,3) -- (13,3); 
          
    \draw[color=gray] (6,2.6) -- (13,2.6); 
    
    \draw[color=yellow] (11,2.8) -- (13,2.8); 
    
        \draw[color=green] (6,1) -- (13,1); 
           
           \draw[color=blue] (11,0.8) -- (13,0.8);

             \draw[color=black] (11,-0.6) -- (13,-0.6); 
             
             \draw[dashed] (1,-2) -- (1,10);
              \draw[dashed] (6,-2) -- (6,10);
               \draw[dashed] (11,-2) -- (11,10);

   \node[draw=red,fill=white,circle,inner sep=1pt] at (1,9.1) {};
   
   \node[draw=red,fill=red,circle,inner sep=1pt] at (1,3) {};
   
   \node[draw=gray,fill=white,circle,inner sep=1pt] at (1,8.5) {};
 
\node[draw=gray,fill=gray,circle,inner sep=1pt] at (1,6.8) {};

  \node[draw=gray,fill=white,circle,inner sep=1pt] at (6,6.8) {};
 
\node[draw=gray,fill=gray,circle,inner sep=1pt] at (6,2.6) {};

  \node[draw=yellow,fill=white,circle,inner sep=1pt] at (1,8.8) {};

  \node[draw=yellow,fill=white,circle,inner sep=1pt] at (1,8.8) {};
 
\node[draw=yellow,fill=yellow,circle,inner sep=1pt] at (1,4.6) {};

  \node[draw=yellow,fill=white,circle,inner sep=1pt] at (11,4.6) {};
 
\node[draw=yellow,fill=yellow,circle,inner sep=1pt] at (11,2.8) {};

  \node[draw=green,fill=white,circle,inner sep=1pt] at (6,7) {};
 
\node[draw=green,fill=green,circle,inner sep=1pt] at (6,1) {};

  \node[draw=black,fill=white,circle,inner sep=1pt] at (6,6.8) {};
 
\node[draw=blue,fill=white,circle,inner sep=1pt] at (6,6.6) {};

 \node[draw=blue,fill=white,circle,inner sep=1pt] at (11,4.8) {};
  
    \node[draw=blue,fill=blue,circle,inner sep=1pt] at (6,4.8) {};

  \node[draw=blue,fill=white,circle,inner sep=1pt] at (11,5) {};
  \node[draw=blue,fill=blue,circle,inner sep=1pt] at (11,0.8) {};
 
\node[draw=black,fill=black,circle,inner sep=1pt] at (11,-0.6) {};

 \draw (-0.2,-2.4)   node {$
 t_{3\nu} = q_\nu^{\{1\}} = $};
 
 \draw(-0.2,-3.2)   node
 {$=q_{2\nu}^{\{1,2\}} = q_{2\nu}^{\{1,3\}}
 $}; 
 
  \draw (6,-2.4)   node {$
  t_{3\nu+1} = q_\nu^{\{2\}} =  $};
 
 \draw(6,-3.2)   node
 {$=q_{2\nu}^{\{2,3\}} = q_{2\nu+1}^{\{1,2\}}
 $};

  \draw (11.4,-2.4)   node {$
   t_{3\nu+2} = q_\nu^{\{3\}} =  $};
 
 \draw(11.4,-3.2)   node
 {$= q_{2\nu+1}^{\{1,3\}} = q_{2\nu+1}^{\{2,3\}}
 $};

    \end{tikzpicture}
    
         \caption{   
       Zum Satz 2 und  seinem Beweis für  $k=3, n = 6$ 
           }

\end{figure}
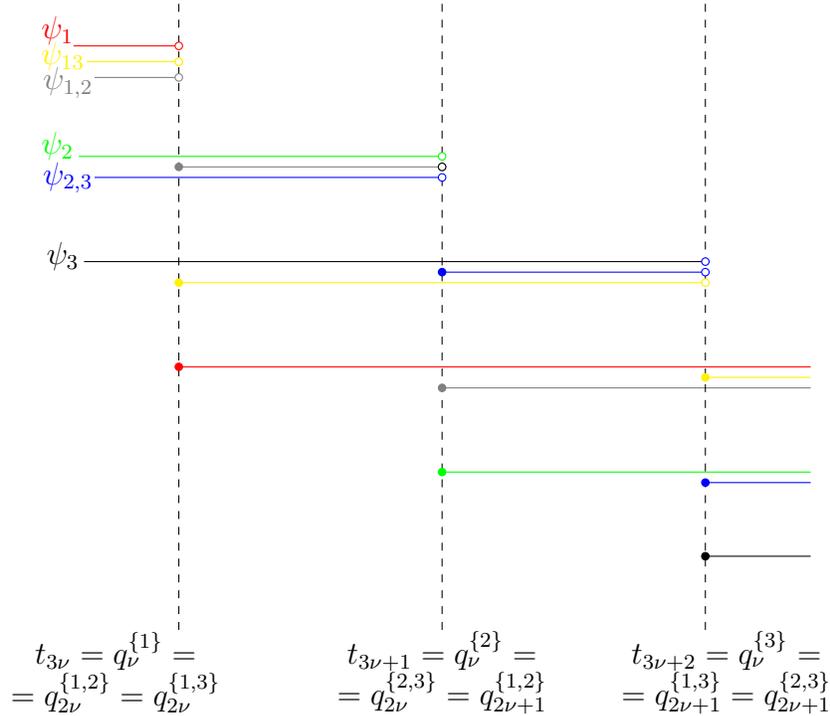
    
  \label{fig1}

\vskip+0.3cm
Wir behaupten, dass   es nur die Permutationen der Form
$$
\psi_{\frak{a}_1^l} (t)>\!\psi_{{\frak{a}_2^l}} (t)>...>\! \psi_{\frak{a}_n^l}(t),\,\,\,\,\, 0\le l \le k-1
$$
für alle hinreichend großen $t$
gibt (die Elemente $\frak{a}_j^l$ sind in (\ref{permuta}) definiert).
   
   Aus der Symmetrie folgt, dass wir uns auf den  folgenden Fall   beschränken können. Nehmen wir an, dass $ t =  t_{k\nu}-1$. Dann müssen wir die Ungleichungen
  $$
  \psi_{\frak{a}_1} (t)>\psi_{\frak{a}_2} (t)>\psi_{\frak{a}_3}(t)>...  >\psi_{\frak{a}_n} (t),\,\,\,\,\, t =  t_{k\nu}-1
  $$ 
   oder
\begin{equation}\label{z1}
\psi_{1} (t)>\psi_{1,k} (t)>...> \psi_{1,2}(t)>\psi_{2}(t)>\psi_{2,k} (t)>... >\psi_{k-1}(t)>\psi_{k-1,k}(t)>\psi_{k} (t),\,\,\, t =  t_{k\nu}-1
\end{equation}
beweisen.   
 Um (\ref{z1}) zu beweisen,  genügt es,   nur die Ungleichungen
\begin{equation}\label{zw3}
 \psi_{i} (t) > \psi_{i,k}(t),\,\,\,\, i= 1,...,k-1,
\end{equation}
\begin{equation}\label{zw31}
  \psi_{i,j}(t)>   \psi_{i,j-1}(t),\,\,\,\, i= 1,...,k-2,\,\, i {+2}  \le j \le k,
\end{equation}
und
\begin{equation}\label{zw32}
 \psi_{i, i+1} (t) > \psi_{i+1}(t),\,\,\,\,\, 1\le i\le k-1
\end{equation}
für $ t  =  t_{k\nu}-1$
zu
beweisen.

Wir beweisen     (\ref{zw3}).
Aus (\ref{zero})  folgen die Gleichungen
 $$
 \psi_i(t_{k\nu}-1) = 
  \psi_i(t_{k(\nu-1)+i-1})=
 \psi_{k}(q^{(i)}_{\nu-1})=
 \frac{1}{q_\nu^{(i)}+q_{\nu-1}^{(i)}/\alpha_{\nu+1}^{(i)}}
 $$
 und
 $$
 \psi_{i,k}(t_{k\nu}-1)
 = 
  \psi_{i,k}(t_{k(\nu-1)+k-1})= \psi_{i,k}(q^{(i,k)}_{2\nu-1})=
   \frac{1}{q_{2\nu}^{(i,k)}+ q_{2\nu-1}^{(i,k)}/\alpha_{2\nu+1}^{(i,k)}}=
  \frac{1}{q_\nu^{(i)}+ q_{2\nu-1}^{(i,k)}/\alpha_{2\nu+1}^{(i,k)}},
 $$
 weil
 $q^{(i)}_{\nu-1}=t_{k(\nu-1)+i-1}$,
 $q_{2\nu-1}^{(i,k)} =  t_{k(\nu-1)+k-1}$
  und $ q^{(i,k)}_{2\nu}= q_\nu^{(i)} =t_{k\nu+i-1}$ gelten
  (wir haben  (\ref{TK}) benutzt).
 Nun gilt
 $$
  \frac{
 \psi_{i}(t_{k\nu}-1)}{ \psi_{i,k}(t_{k\nu}-1)}=
 \frac
 {
 q_\nu^{(i)}+ q_{2\nu-1}^{(i,k)}/\alpha_{2\nu+1}^{(i,k)}
 }
 {
 q_\nu^{(i)}+q_{\nu-1}^{(i)}/\alpha_{\nu+1}^{(i)}
 } >1,
 $$
aufgrund der folgenden asymptotischen Ungleichung:
$$
\frac{
 q_{2\nu-1}^{(i,k)}}{\alpha_{2\nu+1}^{(i,k)}}\sim
  \frac{q_{2\nu-1}^{(i,k)}}{a_{2\nu+1}^{(i,k)}}\sim
  \frac{t_{k(\nu-1)+k-1}t_{k\nu+i-1}}{t_{k\nu+k-1}}
 >
  \frac{t_{k(\nu-1)+i-1}t_{k\nu+i-1}}{t_{k(\nu+1)+i-1}}
 \sim
  \frac{
 q_{\nu-1}^{(i)}}{a_{\nu+1}^{(i)}}
 \sim
 \frac{
 q_{\nu-1}^{(i)}}{\alpha_{\nu+1}^{(i)}},\,\,\,\,\nu\to \infty
$$
 (hier nutzen wir $ a_{2\nu+1}^{(i,k)} \sim \frac{t_{k\nu+k-1}}{t_{k\nu+i-1}}$, 
 $a_{\nu+1}^{(i)} \sim \frac{t_{k(\nu+1)+i-1}}{t_{k\nu+i-1}}$ und (\ref{z00})).
 Die Ungleichung  (\ref{zw3}) ist  damit bewiesen.

 Nun beweisen wir    (\ref{zw31}). Wir haben 
 $$
 \psi_{i,j}(t_{k\nu}-1)
 = 
  \psi_{i,j}(t_{k(\nu-1)+j-1})= \psi_{i,j}(q^{(i,j)}_{2\nu-1})=
   \frac{1}{q_{2\nu}^{(i,j)}+ q_{2\nu-1}^{(i,j)}/\alpha_{2\nu+1}^{(i,j)}},
 $$
 weil die Gleichung
 $q_{2\nu-1}^{(i,j)} =  t_{k(\nu-1)+j-1}$ gilt.
 Wir müssen die Ungleichung
 $$
  \psi_{i,j}(t_{k\nu}-1) =   \frac{1}{q_{2\nu}^{(i,j)}+ q_{2\nu-1}^{(i,j)}/\alpha_{2\nu+1}^{(i,j)}}
  >
    \psi_{i,j-1}(t_{k\nu}-1) =   \frac{1}{q_{2\nu}^{(i,j-1)}+ q_{2\nu-1}^{(i,j-1)}/\alpha_{2\nu+1}^{(i,j-1)}}
    \,\,\,\text{mit}\,\,\, q_{2\nu}^{(i,j)} =q_{2\nu}^{(i,j-1)} = t_{k\nu+i-1}
 $$
 beweisen. Nun gilt
 $
  \frac{
 \psi_{i,j}(t_{k\nu}-1)}{ \psi_{i,j-1}(t_{k\nu}-1)}
  >1
 $,
 weil aus  
 $
  q_{2\nu-1}^{(i,j-1)}= t_{k(\nu-1)+j-2}
,
   q_{2\nu-1}^{(i,j)}= t_{k(\nu-1)+j-1}
 ,$
 $ a_{2\nu+1}^{(i,j)} \sim \frac{t_{k\nu+j-1}}{t_{k\nu+i-1}}$ und (\ref{dop}) die asymptotische Ungleichung
$$
 \frac{
 q_{2\nu-1}^{(i,j-1)}}{\alpha_{2\nu+1}^{(i,j-1)}}\sim
  \frac{q_{2\nu-1}^{(i,j-1)}}{a_{2\nu+1}^{(i,j-1)}}\sim
  \frac{t_{k(\nu-1)+j-2}t_{k\nu+i-1}}{t_{k\nu+j-2}}
 >
  \frac{t_{k(\nu-1)+j-1}t_{k\nu+i-1}}{t_{k\nu+j-1}}
 \sim
  \frac{
 q_{2\nu-1}^{(i,j)}}{a_{2\nu+1}^{(i,j)}}
 \sim
 \frac{
 q_{2\nu-1}^{(i,j)}}{\alpha_{2\nu+1}^{(i,j)}},\,\,\,\,\nu\to \infty
  $$
  folgt.
 Die Ungleichung  (\ref{zw31}) ist  dadurch bewiesen.

Die Ungleichung    (\ref{zw32}) ist klar. Tatsächlich sind
$$ t_{k\nu+i}= q^{(i,i+1)}_{2\nu-1} = q^{(i+1)}_{\nu-1},\,\,\,
{ \psi_{i, i+1} (t_{k\nu}-1) } ={ \psi_{i, i+1} (t_{k\nu+i}) },
\,\,\,
{ \psi_{i+1} (t_{k\nu}-1) } ={ \psi_{ i+1} (t_{k\nu+i}) }
$$
und
$$
\frac{ \psi_{i, i+1} (t_{k\nu}-1)  }{ \psi_{i+1} (t_{k\nu}-1) } =
\frac{ \psi_{i, i+1} (q^{(i,i+1)}_{2\nu-1})  }{ \psi_{i+1} (q^{(i+1)}_{\nu-1})} =
\frac
{q^{(i+1)}_{\nu}(1+{\alpha^{(i+1)*}_{\nu]}}/{\alpha^{(i+1)}_{\nu+1}})}
{q^{(i,i+1)}_{2\nu}(1+{\alpha^{(i,i+1)*}_{2\nu}}/{\alpha^{(i,i+1)}_{2\nu+1}})} \ge
\frac
{q^{(i+1)}_{\nu} }
{2q^{(i,i+1)}_{2\nu} }=
\frac
{t_{k\nu+i}}
{2t_{k\nu+i-1} }\to \infty, \nu\to \infty.
$$
Satz  2 its somit vollständig gezeigt.$\Box$

\end{document}